\documentclass[11pt]{amsart}
\usepackage{times}      
\usepackage{latexsym}   
\usepackage{amssymb}    
\usepackage{amsmath}    
\usepackage{amsthm}     
\usepackage{mathrsfs}   
\newcommand{\bd}{\mathbf{D}}
\newcommand{\bj}{\mathbf{J}}

\newcommand{\id}{\mbox{Id}}

\newcommand{\lp}{\left(}
\newcommand{\rp}{\right)}
\newcommand{\lc}{\left[}
\newcommand{\rc}{\right]}
\newcommand{\lcl}{\left\{}
\newcommand{\rcl}{\right\}}
\newcommand{\lln}{\left|}
\newcommand{\rrn}{\right|}

\newcommand{\eps}{\varepsilon}

\newcommand{\ga}{\gamma}

\newcommand{\oom}{\Omega}

\newcommand{\beq}{\begin{equation}}
\newcommand{\eeq}{\end{equation}}
\newcommand{\bea}{\begin{eqnarray}}
\newcommand{\eea}{\end{eqnarray}}
\newcommand{\beas}{\begin{eqnarray*}}
\newcommand{\eeas}{\end{eqnarray*}}

\def\msh{{\mathscr H}}

\def\cC{{\mathcal C}}

\def\mq{{\mathbb  Q}}
\def\me{{\mathbb  E}}

\def\mr{{\mathbb  R}}

\def\mn{{\mathbb  N}}

\def\mp{{\mathbb  P}}

\newcommand{\cac}{{\mathcal C}}

\newcommand{\cf}{{\mathcal F}}

\newcommand{\ch}{{\mathcal H}}

\newcommand{\crr}{{\mathcal R}}

\newcommand{\D}{{\mathbb D}}
\newcommand{\EE}{{\mathbb E}}

\newcommand{\PP}{{\mathbb P}}

\newcommand{\R}{{\mathbb R}}

\newtheorem{theorem}{Theorem}[section]

\newtheorem{definition}[theorem]{Definition}

\newtheorem{hypothesis}[theorem]{Hypothesis}
\newtheorem{lemma}[theorem]{Lemma}

\newtheorem{proposition}[theorem]{Proposition}
\theoremstyle{remark}
\newtheorem{remark}[theorem]{Remark}
\theoremstyle{remark}
\newtheorem{example}[theorem]{Example}
\theoremstyle{remark}

\newtheorem{assumption}[theorem]{Assumption}
\newtheorem{foo}[theorem]{Remarks}

\def\msh{{\mathscr H}}

\def\cC{{\mathcal C}}

\def\mq{{\mathbb  Q}}
\def\me{{\mathbb  E}}

\def\mr{{\mathbb  R}}

\def\mn{{\mathbb  N}}

\def\mp{{\mathbb  P}}

\numberwithin{equation}{section}

\title[]{On small time asymptotics for rough differential equations driven by fractional Brownian motions}
\author{Fabrice Baudoin}
\address{Department of Mathematics\\
       Purdue University \\ West Lafayette, IN 47907.}
\email{fbaudoin\@@math.purdue.edu}
\thanks{The first author of this research was supported in part by
NSF Grant DMS 0907326}

\author{Cheng Ouyang}
\address{Dept. Mathematics, Statistics and Computer Science\\ University of Illinois at Chicago\\ Chicago, IL 60607.}
\email{couyang@math.uic.edu}

\begin{document}

\begin{abstract}

We survey existing results concerning the study in small times of the density of the solution of a rough differential equation driven by fractional Brownian motions. We also slightly improve existing results and discuss some possible applications to mathematical finance.
\end{abstract}

\maketitle

\begin{center}
\textit{In memory of Peter Laurence}
\end{center}

\tableofcontents

\section{Introduction} 

In this paper, our main goal is to survey some existing results concerning the small-time asymptotics of the density of rough differential equations driven by fractional Brownian motions. Even though we do not  claim any new results, we slightly  improve some of the existing ones and also point out  some possible connections to finance. We also hope, it will be useful for the reader to have, in one place, the most recent results concerning the small-time asymptotics questions related to rough differential equations driven by fractional Brownian motions.
Our discussion will mainly be based on one hand on the papers \cite{BO1,BO2,BOZ1} by the two present authors and on the other hand on the papers \cite{Inahama1,Inahama2} by Inahama.

\

Random dynamical systems are a well established modeling tool for a variety of natural phenomena ranging from physics (fundamental and phenomenological) to chemistry and more recently to biology, economy, engineering sciences and mathematical finance. In many interesting models the lack of any regularity of the external inputs of the differential equation as functions of time is a technical difficulty that hampers their mathematical analysis. The theory of rough paths has been initially developed by T. Lyons \cite{Ly} in the 1990's to provide a framework to analyze a large class of driven differential equations and the precise relations between the driving signal and the output (that is the state, as function of time, of the controlled system). 

Rough paths theory provides a perfect framework to study differential equations driven by Gaussian processes (see \cite{FV1}). In particular, using rough paths theory, we may define solutions of stochastic differential equations driven by a fractional Brownian motion with a parameter $H>1/4$ (see \cite{CQ}). Let us then consider the equation
\begin{align}\label{SDEgaussianintro}
X_t^{x}=x+\int_0^t V_0(X^x_s)ds+\sum_{i=1}^d \int_0^t V_i(X_s^x)dB^i_s,
\end{align}
where $x \in \mathbb{R}^n$, $V_0,V_1,\cdots,V_d$ are bounded smooth vector fields and $(B_t)_{t\ge 0}$ is a $d$-dimensional  fractional Brownian motion with Hurst parameter $H\in(\frac{1}{4},1)$. A first basic question is the existence of a smooth density with respect to the Lebesgue measure for the random variable $X_t^{x}$, $t>0$.  After multiple works, it is now understood that the answer to this question is essentially the same as the one for stochastic differential equations driven by Brownian motions: the random variable  $X_t^{x}$ admits a smooth density with respect to the Lebesgue measure if H\"ormander's condition is satisfied at $x$ .  More precisely, if $I=(i_1,\ldots,i_k) \in \{ 0,\ldots, d \}^k$, we denote by $V_I$ the Lie commutator defined by
\[
V_I = [V_{i_1},[V_{i_2},\ldots,[V_{i_{k-1}}, V_{i_{k}}]\ldots],
\]
and 
\[
d(I)=k+n(I),
\]
where $n(I)$ is the number of 0 in the word $I$. The basic and fundamental result concerning the existence of a density for stochastic differential equations driven by fractional Brownian motions is the following:

\begin{theorem}[\cite{baudoin14,CF,CLL,H-P}]\label{theo:main}
Assume $H>\frac{1}{4}$ and assume that, at some $x \in \mathbb{R}^n$, there exists $N$ such
that
\begin{equation}
\label{Hormander assumption} \mathbf{span} \{ V_I(x), d(I)\le N \}=\mathbb{R}^n\;.
\end{equation}
Then, for any $t>0$, the law of the random variable $X_t^{x}$ has a smooth density $p_t(x,y)$ with respect to the Lebesgue measure on $\mathbb{R}^n$.
\end{theorem}

Once the existence and smoothness of the density is established, it is natural to study  properties of this density. In particular, we are interested here in small-time asymptotics, that is the analysis of $p_t(x,y)$ when $ t \to 0$. Based on the results in the Brownian motion case \cite{azencott,benarous}, and taking into account the scaling property of the fractional Brownian motion, the following expansion is somehow expected when $x,y$ are not in the cut-locus one of each other:
\begin{align}\label{kjh}
p_t(x,y)=\frac{1}{(t^H)^{d}}e^{-\frac{d^2(x,y)}{2t^{2H}}}\bigg(\sum_{i=0}^N c_i(x,y)t^{2iH}+r_{N+1}(t,x,y)t^{2(N+1)H}\bigg).
\end{align}

Our goal is to discuss here the various assumptions under which such expansion is known to be true and also discuss possible variations. The approach to study the problem is similar to the case of Brownian motion, the main difficulty to overcome is to study the Laplace method on the path space of the fractional Brownian motion (see \cite{benarous2} for the Brownian case).

The paper is organized as follows. In Section 2 we give some basic results of the theory of rough paths and of the Malliavin calculus tools that will be needed.  In Section 3, we prove a Varadhan's type small time asymptotics for $\ln p_t(x,y)$. The discussion is mainly based on \cite{BOZ1}. In Section 4, we study sufficient conditions under which the above expansion \eqref{kjh}  is valid. Our discussion is based on \cite{BO1,Inahama1,Inahama2}. Finally, in Section 5, we discuss some models in mathematical finance where the asymptotics of the density for rough differential equations may play an important role.

\section{Preliminary material}

For some fixed $H>\frac{1}{4}$, we consider $(\oom,\cf,\PP)$ the canonical probability space associated with the fractional
Brownian motion (in short fBm) with Hurst parameter $H$. That is,  $\oom=\cac_0([0,1])$ is the Banach space of continuous functions
vanishing at zero equipped with the supremum norm, $\cf$ is the Borel sigma-algebra and $\PP$ is the unique probability
measure on $\oom$ such that the canonical process $B=\{B_t=(B^1_t,\ldots,B^d_t), \; t\in [0,1]\}$ is a fractional Brownian motion with Hurst
parameter $H$.
In this context, let us recall that $B$ is a $d$-dimensional centered Gaussian process, whose covariance structure is induced by
\begin{align}\label{covariance}
R\left( t,s\right) :=\EE\lc  B_s^j \, B_t^j\rc
=\frac{1}{2}\left( s^{2H}+t^{2H}-|t-s|^{2H}\right),
\quad
s,t\in[0,1] \mbox{ and } j=1,\ldots,d.
\end{align}
In particular it can be shown, by a standard application of Kolmogorov's criterion, that $B$ admits a continuous version
whose paths are $\ga$-H\"older continuous for any $\ga<H$.

\subsection{Rough paths theory}

In this section, we recall some basic results in rough paths theory. More details can be found in the monographs \cite{FV-bk} and \cite{LQ}.
For $N\in\mathbb{N}$, recall that the truncated algebra $T^{N}(\mathbb{R}%
^{d})$ is defined by
$$
T^{N}(\mathbb{R}^{d})=\bigoplus_{m=0}^{N}(\mathbb{R}%
^{d})^{\otimes m},
$$
with the convention $(\mathbb{R}^{d})^{\otimes
0}=\mathbb{R}$. The set $T^{N}(\mathbb{R}^{d})$ is equipped with a straightforward
vector space structure plus an multiplication $\otimes$. 
Let $\pi_{m}$ be the projection on the $m$-th tensor level. Then
$(T^{N}(\mathbb{R}^{d}),+,\otimes)$ is an associative algebra with unit
element $\mathbf{1} \in(\mathbb{R}^{d})^{\otimes0}$.

\smallskip

For $s<t$ and $m\geq2$, consider the simplex $\Delta_{st}^{m}=\{(u_{1}%
,\ldots,u_{m})\in\lbrack s,t]^{m};\,u_{1}<\cdots<u_{m}\} $, while the
simplices over $[0,1]$ will be denoted by $\Delta^{m}$. A continuous map
$\mathbf{x}:\Delta^{2}\rightarrow T^{N}(\mathbb{R}^{d})$ is called a
multiplicative functional if for $s<u<t$ one has $\mathbf{x}_{s,t}%
=\mathbf{x}_{s,u}\otimes\mathbf{x}_{u,t}$. An important example arises from
considering paths $x$ with finite variation: for $0<s<t$ we set
\begin{equation}
\mathbf{x}_{s,t}^{m}=\sum_{1\leq i_{1},\ldots,i_{m}\leq d}\biggl( \int%
_{\Delta_{st}^{m}}dx^{i_{1}}\cdots dx^{i_{m}}\biggr) \,e_{i_{1}}\otimes
\cdots\otimes e_{i_{m}}, \label{eq:def-iterated-intg}%
\end{equation}
where $\{e_{1},\ldots,e_{d}\}$ denotes the canonical basis of $\mathbb{R}^{d}%
$, and then define the truncated \textit{signature} of $x$ as
\[
S_{N}(x):\Delta^{2}\rightarrow T^{N}(\mathbb{R}^{d}),\qquad(s,t)\mapsto
S_{N}(x)_{s,t}:=1+\sum_{m=1}^{N}\mathbf{x}_{s,t}^{m}.
\]
The function $S_{N}(x)$ for a smooth function $x$ will be our typical example of multiplicative functional. Let us stress the fact that those elements take values in the strict subset
$G^{N}(\mathbb{R}^{d})\subset T^{N}(\mathbb{R}^{d})$,  called free nilpotent group of step $N$, and is equipped with the classical Carnot-Caratheodory norm which we simply denote by $|\cdot|$. For a path $\mathbf{x}\in\cC([0,1],G^{N}(\R^d))$, the $p$-variation norm of $\mathbf{x}$ is defined to be
\begin{align*}
\|\mathbf{x}\|_{p-{\rm var}; [0,1]}=\sup_{\Pi \subset [0,1]}\left(\sum_i |\mathbf{x}_{t_i}^{-1}\otimes \mathbf{x}_{t_{i+1}}|^p\right)^{1/p}
\end{align*}
where the supremum is taken over all subdivisions $\Pi$ of $[0,1]$.

\smallskip

With these notions in hand, let us briefly define what we mean by geometric rough path (we refer to \cite{FV-bk,LQ} for a complete overview): for $p\geq 1$, an element $x: [0,1]\to G^{\lfloor p \rfloor}(\mr^d)$ is said to be a geometric rough path if it is the $p$-var limit of a sequence $S_{\lfloor p \rfloor}(x^{m})$.  In particular, it is an element of the space
$$\cC^{p-{\rm var}; [0,1]}([0,1], G^{\lfloor p \rfloor}(\R^d))=\{\mathbf{x}\in \cC([0,1], G^{\lfloor p \rfloor}(\R^d)): \|\mathbf{x}\|_{p-{\rm var}; [0,1]}<\infty\}.
$$

Let $\mathbf{x}$ be a geometric $p$-rough path with its approximating sequence $x^m$, that is, $x^m$ is a sequence of smooth functions such that $\mathbf{x}^m=S_{\lfloor p\rfloor}(x^m)$ converges to $\mathbf{x}$ in the $p$-var norm. Fix any $1\leq q\leq p$ so that $p^{-1}+ q^{-1}>1$ and pick any $h\in\cC^{q-{\rm var}}([0,1], \R^d)$. One can define the translation of $\mathbf{x}$ by h, denoted by $T_h(\mathbf{x})$ by
$$T_h(\mathbf{x})=\lim_{n\to\infty}S_{\lfloor p\rfloor}({x}^m+h).$$
It can be shown that $T_h(\mathbf{x})$ is an element in $\cC^{p-{\rm var}}([0,1], G^{\lfloor p \rfloor}(\R^d))$. Moreover, one can show that $T_h(\mathbf{x})$ uniformly continuous in $h$ and $\mathbf{x}$ on bounded sets.


\begin{remark}
A typical situation of the above translation of $\mathbf{x}$ by $h$ in the present paper is when $\mathbf{x}=\mathbf{B}$, the fractional Brownian motion lifted as a rough path, and $h$ is a Cameron-Martin element of $B$. In this case, we simply denote $T_h(\mathbf{B})=B+h.$ 
\end{remark}

According to the considerations above, in order to prove that a lift of a $d$-dimensional fBm as a geometric rough path exists it is sufficient to build enough iterated integrals of $B$ by a limiting procedure. Towards this aim, a lot of the information concerning $B$ is encoded in the rectangular increments of the covariance function $R$ (defined by \eqref{covariance}), which are given by
\begin{equation*}
R_{uv}^{st} \equiv \EE\lc (B_t^1-B_s^1) \, (B_v^1-B_u^1) \rc.
\end{equation*}
We then call 2-dimensional $\rho$-variation of $R$ the quantity
\begin{equation*}
V_{\rho}(R)^{\rho} \equiv
\sup \lcl
\lp \sum_{i,j} \lln R_{s_{i} s_{i+1}}^{t_{j}t_{j+1}} \rrn^{\rho} \rp^{1/\rho}; \, (s_i), (t_j)\in \Pi
\rcl,
\end{equation*}
where $\Pi$ stands again for the set of partitions of $[0,1]$. The following result is now well known for fractional Brownian motion \cite{CQ,FV1}:
\begin{proposition}\label{prop:fbm-rough-path}
For a fractional Brownian motion with Hurst parameter $H$, we have $V_{\rho}(R)<\infty$ for all $\rho\ge1/(2H)$. Consequently, for $H>1/4$ the process $B$ admits a lift $\mathbf{B}$ as a geometric rough path of order $p$ for any $p>1/H$.
\end{proposition}



\subsection{Malliavin Calculus} We introduce the basic framework of Malliavin calculus in this subsection.  The reader is invited to consult the corresponding chapters in  \cite{Nu06} for further details. Let $\mathcal{E}$ be the space of $\mathbb{R}^d$-valued step
functions on $[0,1]$, and $\mathcal{H}$  the closure of
$\mathcal{E}$ for the scalar product:
\[
\langle (\mathbf{1}_{[0,t_1]} , \cdots ,
\mathbf{1}_{[0,t_d]}),(\mathbf{1}_{[0,s_1]} , \cdots ,
\mathbf{1}_{[0,s_d]}) \rangle_{\mathcal{H}}=\sum_{i=1}^d
R(t_i,s_i).
\]
We denote by $K^*_H$ the isometry between $\mathcal{H}$ and $L^2([0,1])$.
When $H>\frac{1}{2}$ it can be shown that $\mathbf{L}^{1/H} ([0,1], \mathbb{R}^d)
\subset \mathcal{H}$, and when $\frac{1}{4}<H<\frac{1}{2}$ one has
$$C^\gamma\subset \mathcal{H}\subset L^2([0,1])$$
for all $\gamma>\frac{1}{2}-H$.

We remark that $\ch$ is the reproducing kernel Hilbert space for $B$. Let $\msh_H$ be the Cameron-Martin space of $B$,
one proves that the operator $\crr:=\crr_H :\ch \rightarrow \msh_H$ given by
\begin{equation}\label{eq:def-R}
\crr \psi := \int_0^\cdot K_H(\cdot,s) [K^*_H \psi](s)\, ds
\end{equation}
defines an isometry between $\ch$ and $\msh_H$. 
Let us now quote from \cite[Chapter 15]{FV-bk} a result relating the 2-d regularity of $R$ and the regularity of $\msh_H$.
\begin{proposition}\label{prop:imbed-bar-H}
Let $B$ be a fBm with Hurst parameter $\frac{1}{4}<H<\frac{1}{2}$. Then one has $\msh_H\subset \cac^{\rho-{\rm var}}$ for $\rho>(H+1/2)^{-1}$. Furthermore, the following quantitative bound holds:
\begin{equation*}
\Vert h \Vert_{\msh_H} \ge \frac{\Vert h \Vert_{\rho-{\rm var}}}{(V_\rho(R))^{1/2}}.
\end{equation*}
\end{proposition}
\begin{remark}
The above proposition shows that for fBm we have $\msh_H\subset \cac^{\rho-{\rm var}}$ for $\rho>(H+1/2)^{-1}$. Hence an integral of the form $\int h \, dB$ can be interpreted in the Young sense by means of $p$-variation techniques.
\end{remark}

A $\mathcal{F}$-measurable real
valued random variable $F$ is then said to be cylindrical if it can be
written, for a given $n\ge 1$, as
\begin{equation*}
F=f\lp  B(\phi^1),\ldots,B(\phi^n)\rp=
f \Bigl( \int_0^{1} \langle \phi^1_s, dB_s \rangle ,\ldots,\int_0^{1}
\langle \phi^n_s, dB_s \rangle \Bigr)\;,
\end{equation*}
where $\phi^i \in \mathcal{H}$ and $f:\mathbb{R}^n \rightarrow
\mathbb{R}$ is a $C^{\infty}$ bounded function with bounded derivatives. The set of
cylindrical random variables is denoted $\mathcal{S}$.

The Malliavin derivative is defined as follows: for $F \in \mathcal{S}$, the derivative of $F$ is the $\mathbb{R}^d$ valued
stochastic process $(\mathbf{D}_t F )_{0 \leq t \leq 1}$ given by
\[
\mathbf{D}_t F=\sum_{i=1}^{n} \phi^i (t) \frac{\partial f}{\partial
x_i} \left( B(\phi^1),\ldots,B(\phi^n)  \right).
\]
More generally, we can introduce iterated derivatives. If $F \in
\mathcal{S}$, we set
\[
\mathbf{D}^k_{t_1,\ldots,t_k} F = \mathbf{D}_{t_1}
\ldots\mathbf{D}_{t_k} F.
\]
For any $p \geq 1$, it can be checked that the operator $\mathbf{D}^k$ is closable from
$\mathcal{S}$ into $\mathbf{L}^p(\oom;\mathcal{H}^{\otimes k})$. We denote by
$\mathbb{D}^{k,p}$ the closure of the class of
cylindrical random variables with respect to the norm
\[
\left\| F\right\| _{k,p}=\left( \mathbb{E}\left( F^{p}\right)
+\sum_{j=1}^k \mathbb{E}\left( \left\| \mathbf{D}^j F\right\|
_{\mathcal{H}^{\otimes j}}^{p}\right) \right) ^{\frac{1}{p}},
\]
and
\[
\mathbb{D}^{\infty}=\bigcap_{p \geq 1} \bigcap_{k
\geq 1} \mathbb{D}^{k,p}.
\]

\begin{definition}\label{non-deg}
Let $F=(F^1,\ldots , F^n)$ be a random vector whose components are in $\mathbb{D}^\infty$. Define the Malliavin matrix of $F$ by
$$\gamma_F=(\langle \mathbf{D}F^i, \mathbf{D}F^j\rangle_{\ch})_{1\leq i,j\leq n}.$$
Then $F$ is called  {\it non-degenerate} if $\gamma_F$ is invertible $a.s.$ and
$$(\det \gamma_F)^{-1}\in \cap_{p\geq1}L^p(\Omega).$$
\end{definition}
\noindent
It is a classical result that the law of a non-degenerate random vector $F=(F^1, \ldots , F^n)$ admits a smooth density with respect to the Lebesgue measure on $\mr^n$. Furthermore, the following integration by parts formula allows to get more quantitative estimates:

\begin{proposition}\label{th: IBP}
Let $F=(F^1,...,F^n)$ be a non-degenerate random vector whose components are in $\D^\infty$, and $\gamma_F$ the Malliavin matrix of $F$. Let $G\in\D^\infty$ and $\varphi$ be a function in the space $C_p^\infty(\mr^n)$. Then for any multi-index $\alpha\in\{1,2,...,n\}^k, k\geq 1$, there exists an element $H_\alpha\in\D^\infty$ such that
$$\me[\partial_\alpha \varphi(F)G]=\me[\varphi(F)H_\alpha].$$
Moreover, the elements $H_\alpha$ are recursively given by
\begin{align*}
&H_{(i)}=\sum_{j=1}^{d}\delta\left(G(\gamma_F^{-1})^{ij}\mathbf{D}F^j\right)\\
&H_\alpha=H_{(\alpha_k)}(H_{(\alpha_1,..., \alpha_{k-1})}),
\end{align*}
and for $1\leq p<q<\infty$ we have
$$\|H_\alpha\|_{L^p}\leq C_{p,q}\|\gamma_F^{-1}\mathbf{D}F\|^k_{k, 2^{k-1}r}\|G\|_{k,q},$$
where $\frac{1}{p}=\frac{1}{q}+\frac{1}{r}$.
\end{proposition}
\begin{remark}\label{est H}
By the estimates for $H_\alpha$ above, one can conclude that there exist constants $\beta, \gamma>1$ and integers $m, r$ such that
\begin{align*}
\|H_\alpha\|_{L^p}\leq C_{p,q}\|\det\gamma^{-1}_F\|^m_{L^\beta}\|\mathbf{D}F\|_{k,\gamma}^r\|G\|_{k,q}.
\end{align*}
\end{remark}
\begin{remark}
In what follows, we use $H_\alpha(F,G)$ to emphasize its dependence on $F$ and $G$.
\end{remark}


\subsection{Differential equations driven by fractional Brownian motions}
Let $B$ be a d-dimensional fractional Brownian motion with Hurst parameter $H>\frac{1}{4}$. Fix a small parameter $\eps\in(0,1]$, and consider the solution $X_t^\eps$ to the stochastic differential equation
\begin{align}\label{equ: SDE}
X_t^\eps=x+\eps\sum_{i=1}^d\int_0^tV_i(X_s^\eps)dB_s^i+\int_0^tV_0(\eps, X_s^\eps)ds,
\end{align}
where the vector fields $V_1,\ldots,V_d$ are $C^\infty$-bounded vector fields on $\R^n$ and $V_0(\eps, \cdot)$ is $C^\infty$-bounded uniform in $\eps\in [0,1]$. 


Proposition \ref{prop:fbm-rough-path} ensures the existence of a lift of $B$ as a geometrical rough path. The general rough paths theory (see e.g.~\cite{FV-bk,Gu}) allows thus to state  the following proposition:

\begin{proposition}\label{prop:moments-sdes-rough}
Consider equation (\ref{equ: SDE}) driven by a $d$-dimensional fBm $B$ with Hurst parameter $H>\frac{1}{4}$, and assume that the vector fields $V_i$s are $C^\infty$-bounded. Then

\smallskip

\noindent
\emph{(i)}
For each $\eps\in(0,1]$, equation (\ref{equ: SDE}) admits a unique finite $p$-var continuous solution $X^\eps$ in the rough paths sense, for any $p> \frac{1}{H}$.

\smallskip

\noindent
\emph{(ii)}
For any $\lambda>0$ and $\delta<\frac{1}{p}$ we have
\begin{equation}\label{eq:exp-delta-moments}
\me\left[\exp\lambda\left(\sup_{t\in [0,1], \epsilon \in (0,1]}|X^\eps_t|^\delta\right)\right]<\infty.
\end{equation}
\end{proposition}

Once equation (\ref{equ: SDE}) is solved, the vector $X_t^\eps$ is a typical example of  random variable which can be differentiated in the Malliavin sense. We shall express this Malliavin derivative in terms of the Jacobian $\bj^\eps$ of the equation, which is defined by the relation $$\bj_{t}^{\eps, ij}=\partial_{x_j}X_t^{\eps,i}.$$ Setting $DV_{j}$ for the Jacobian of $V_{j}$ seen as a function from $\R^{n}$ to $\R^{n}$, let us recall that $\bj^\eps$ is the unique solution to the linear equation
\begin{equation}\label{eq:jacobian}
\bj_{t}^\eps = \id_{n} +
\eps\sum_{j=1}^d \int_0^t DV_j (X^{\eps}_s) \, \bj_{s}^\eps \, dB^j_s,
\end{equation}
and that the following results hold true (see \cite{CF} and \cite{NS}  for further details):
\begin{proposition}\label{prop:deriv-sde}
Let $X^\eps$ be the solution to equation (\ref{equ: SDE}) and suppose the $V_i$'s are $C^\infty$-bounded. Then
for every $i=1,\ldots,n$, $t>0$, and $x \in \mathbb{R}^n$, we have $X_t^{\eps,i} \in
\mathbb{D}^{\infty}$ and
\begin{equation*}
\mathbf{D}^j_s X_t^{\eps}= \mathbf{J}^\eps_{st} V_j (X^\eps_s) , \quad j=1,\ldots,d, \quad
0\leq s \leq t,
\end{equation*}
where $\mathbf{D}^j_s X^{\eps,i}_t $ is the $j$-th component of
$\mathbf{D}_s X^{\eps,i}_t$, $\mathbf{J}_{t}^\eps=\partial_{x} X^\eps_t$ and $\bj_{st}^\eps=\bj_{t}^\eps(\bj_{s}^\eps)^{-1}$.

\end{proposition}

\smallskip

Let us now quote the recent result \cite{CLL}, which gives a useful estimate for moments of the Jacobian of rough differential equations driven by Gaussian processes.
\begin{proposition}\label{prop:moments-jacobian}
Consider a  fractional Brownian motion $B$ with Hurst parameter $H>\frac{1}{4}$ and $p>\frac{1}{H}$. Then for any $\eta\ge 1$, there exists a finite constant $c_\eta$ such that the Jacobian $\bj^\eps$ defined at Proposition \ref{prop:deriv-sde} satisfies:
\begin{equation}\label{eq:moments-J-pvar}
\EE\lc  \sup_{\eps\in[0,1]}\Vert \bj^\eps \Vert^{\eta}_{p-{\rm var}; [0,1]} \rc = c_\eta.
\end{equation}
\end{proposition}

Finally, we close the discussion of this section by the following large deviation principle that will be needed later.
Let $\Phi: \msh_H\to \mathcal{C}([0,1],\mathbb{R}^{n})$ be given by solving the ordinary diferential equation
\begin{align}\label{phi}\Phi_t(h)=x+\sum_{i=1}^d\int_0^tV_i(\Phi_s(h))dh_s^i+\int_0^t V_0(0,\Phi_s(h))ds.\end{align}

\begin{theorem}\label{th: LDP}
Let $\Phi$ be given in (\ref{phi}), which is a differentiable mapping from $\msh_H$ to $\mathcal{C}([0,1],\mathbb{R}^{n})$. Denote by $\gamma_{\Phi_1(h)}$ the deterministic Malliavin matrix
 of $\Phi_1(h)$, i.e., $\gamma^{ij}_{\Phi_1(h)}=\langle \bd \Phi_1^i(h), \bd\Phi_1^j(h)\rangle_\ch$, and introduce the following functions on $\mr^n$ and $\mr^n\times\mr$, respectively
$$I(y)=\inf_{\Phi_1(h)=y}\frac{1}{2}\|h\|_{\msh_H}^2, \quad\textnormal{and}\ \  I_R(y,a)=\inf_{\Phi_1(h)=y, \gamma_{\Phi_1(h)}=a}\frac{1}{2}\|h\|_{\msh_H}^2.$$
Recall that $X_1^\eps$ is the solution to equation (\ref{equ: SDE}) and $\gamma_{X_1^\eps}$ is the Malliavin matrix of $X_1^\eps$. Then

(1) $X_1^\eps$ satisfies a large deviation principle with rate function $I(y)$.

(2)The couple $(X_1^\eps,\gamma_{X_1^\eps})$ satisfies a large deviation principle with rate function $I_R(y,a)$.

\end{theorem}



\section{Varadhan asymptotics}

In this section, we are interested in a family of stochastic differential equations driven by fractional Brownian motions $B$ (with Hurst parameter $H>\frac{1}{4}$) of the following form
\begin{align*}
X_t^\eps=x+\eps\sum_{i=1}^d\int_0^tV_i(X_s^\eps)dB_s^i.
\end{align*}
We define a map $\Phi: \msh_H\to \mathcal{C}[0,1]$  by solving the ordinary differential equation
$$\Phi_t(h)=x+\sum_{i=1}^d\int_0^tV_i(\Phi_s(h))dh_s^i.$$
Clearly, we have $X_t^\eps=\Phi_t(\eps \mathbf{B})$. Introduce the following functions on $\mr^n$, which depends on $\Phi$
$$d^2(y)=I(y)=\inf_{\Phi_1(h)=y}\frac{1}{2}\|h\|_{\msh_H}^2,\quad\mathrm{and}\quad d^2_R(y)=\inf_{\Phi_1(h)=y, \det\gamma_{\Phi_1(h)}>0}\frac{1}{2}{\|h\|^2_{\msh_H}}.$$

\bigskip
Throughout the section, we assume that the following assumption Hypothesis \ref{UH condition}  is satisfied. Let us first introduce some notations. Let $\mathcal{A}=\{\emptyset\}\cup\bigcup^{\infty}_{k=1}\{1,2,\cdots,n\}^{k}$ and
$\mathcal{A}_{1}=A\setminus\{\emptyset\}$. We say that $I\in \mathcal{A}$ is a word of length $k$ if $I=(i_1,\cdots,i_k)$
and we write $|I|=k$. If $I=\emptyset$, then we denote $|I|=0$. For any integer $l\ge 1$, we denote by $\mathcal{A}(l)$
the set $\{I\in \mathcal{A}; |I|\le l\}$ and by $\mathcal{A}_{1}(l)$ the set $\{I\in \mathcal{A}_{1}; |I|\le l\}$ .
We also define an operation $\ast$ on $\mathcal{A}$ by
$I\ast J=(i_1,\cdots,i_k,j_1,\cdots,j_l)$ for $I=(i_1,\cdots,i_k)$ and $J=(j_1,\cdots, j_l)$ in $\mathcal{A}$.
We define vector fields $V_{[I]}$ inductively by
\[
V_{[j]}=V_{j}, \quad V_{[I\ast j]}=[V_{[I]}, V_{j}], \quad j=1,\cdots,d
\]

\begin{hypothesis}\label{UH condition}(Uniform hypoelliptic condition)
The vector fields $V_{1},\cdots,V_{d}$ are in $C^{\infty}_{b}(\mathbb{R}^{n})$ and they form a uniform hypoelliptic system 
in the sense that there exist an integer $l$ and a constant $\lambda>0$ such that
\begin{align}\label{UH_condition}
\sum_{I\in \mathcal{A}_{1}(l)}\langle V_{[I]}(x), u\rangle^{2}_{\mathbb{R}^{n}}\ge \lambda \|u\|^{2}
\end{align}
holds for any $x,u\in \mathbb{R}^{n}$
\end{hypothesis}

Under this assumption the main result proved in \cite{BOZ1} is the following Varadhan's type estimate:

\begin{theorem}\label{th: main result}
Let us denote by $p_\eps(y)$ the density of $X_1^\eps$. Then
\begin{align}\label{main claim 1}
\liminf_{\eps\downarrow0}\eps^2\log p_\eps(y)\geq -d^2_R(y),
\end{align}
and
\begin{align}\label{main claim 2}
\limsup_{\eps\downarrow0}\eps^2\log p_\eps(y)\leq -d^2(y).
\end{align}
Moreover, if $$\inf_{\Phi_1(h)=y, \det\gamma_{\Phi_1(y)}>0}\det\gamma_{\Phi_1(h)}>0,$$ then
\begin{align}\label{main claim 3}
\lim_{\eps\downarrow0}\eps^2\log p_\eps(y)=-d^2_R(y).
\end{align}
\end{theorem}

The two key ingredients in proving Theorem \ref{th: main result} are an estimate for the Malliavin derivative $\bd X_1^\eps$ and an estimate of the Malliavin matrix $\gamma_{X_1^\eps}$ of $X_1^\eps$. Building on previous results from \cite{BOZ2},  the following estimates were obtained in \cite{BOZ1} :

\begin{lemma}\label{th: Malliavin est subelliptic}
Assume Hypothesis \ref{UH condition}.  For $H>\frac{1}{4}$, we have
\begin{itemize}
\item[(1)] $\sup_{\eps\in(0,1]}\|X_1^\eps\|_{k,r}<\infty$ for each $k\geq 1$ and $r\geq 1$.
\item[(2)] $\|\gamma_{X_1^\eps}^{-1}\|_r\leq c_r \eps^{-2l}$ for any $r\geq 1$.
\end{itemize}
\end{lemma}

\noindent {\bf Sketch of the proof of (\ref{main claim 1})} Fix $y\in\mr^n$. We only need to show for  $d^2_R(y)<\infty$, since if $d^2_R(y)=\infty$ the statement is trivial. Fix any $\eta>0$ and let $h\in\msh_H$ be such that $\Phi_1(h)=y, \det_{\gamma_{\Phi}}(h)>0$, and $\|h\|^2_{\msh_H}\leq d^2_R(y)+\eta$. Let $f\in C_0^\infty(\mr^n).$ By the Cameron-Martin theorem for fractional Brownian motions, we have
$$\me f(X^\eps_1)=e^{-\frac{\|h\|_{\msh_H}^2}{2\eps^2}}\me f(\Phi_1(\eps B+h))e^\frac{B(h)}{\eps}.$$
Consider then a function $\chi\in C^\infty(\mr), 0\leq \chi\leq 1,$ such that $\chi(t)=0$ if $t\not\in[-2\eta, 2\eta]$, and $\chi(t)=1$ if $t\in[-\eta,\eta]$. Then, if $f\geq 0$, we have
$$\me f(X^\eps_1)\geq e^{-\frac{\|h\|_{\msh_H}+4\eta}{2\eps^2}}\me \chi(\eps B(h))f(\Phi_1(\eps B+h)).$$
Hence, we obtain
\begin{align}\label{lower bound claim1}\eps^2\log p_\eps(y)\geq -(\frac{1}{2}\|h\|_{\msh_H}^2+2\eta)+\eps^2\log\me\big(\chi(\eps B(h))\delta_y(\Phi_1(\eps B+h))\big).\end{align}

On the other hand, we have
$$\me\big(\chi(\eps B(h))\delta_y(\Phi_1(\eps B+h))\big)=\eps^{-n}\me\left(\chi(\eps B(h))\delta_0\left(\frac{\Phi_1(\eps B+h)-\Phi_1(h)}{\eps}\right)\right).$$
Note that
$$Z_1(h)=\lim_{\eps\downarrow 0}\frac{\Phi_1(\eps B+h)-\Phi_1(h)}{\eps}$$
is a $n$-dimensional random vector in the first Wiener chaos with variance $\gamma_{\Phi_1}(h)>0$. Hence $Z_1(h)$ is non-degenerate and  we can then prove that we  obtain
$$\lim_{\eps\downarrow0}\me\left(\chi(\eps B(h))\delta_0\left(\frac{\Phi_1(\eps B+h)-\Phi_1(h)}{\eps}\right)\right)=\me\delta_0(Z_1(h)).$$
Therefore,
$$\lim_{\eps\downarrow0}\eps^2\log\me\big(\chi(\eps B(h))\delta_y(\Phi_1(\eps B+h))\big)=0.$$
Letting $\eps\downarrow0$ in (\ref{lower bound claim1}) we obtain
$$\liminf_{\eps\downarrow0}\eps^2\log p_{\eps}(y)\geq-(\frac{1}{2}\|h\|^2_{\msh_H}+2\eta)\geq -(d^2_R(y)+3\eta).$$
Since $\eta>0$ is arbitrary, this completes the proof. \hfill$\Box$

\noindent{\bf Sketch of the proof of (\ref{main claim 2})}. Fix a point $y\in\mr^n$ and consider a function $\chi\in C_0^\infty(\mr^n), 0\leq\chi\leq1$ such that $\chi$ is equal to one in a neighborhood of $y$. The density of $X_1^\eps$ at point $y$ is given by
$$p_\eps(y)=\me (\chi(X_1^\eps)\delta_y(X_1^\eps)).$$
By the integration by parts formula of Proposition \ref{th: IBP}, we can write
\begin{align*}
\me\chi(X_1^\eps)\delta_y(X_1^\eps)=&\me\left(\mathbf{1}_{\{X_1^\eps>y\}}H_{(1,2,...,n)}(X_1^\eps,\chi(X_1^\eps))\right)\\
\leq&\me|H_{(1,2,...,n)}(X_1^\eps,\chi(X_1^\eps))|\\
=&\me\big(|H_{(1,2,...,n)}(X_1^\eps,\chi(X_1^\eps))|\mathbf{1}_{\{X_1^\eps\in \mathrm{supp}\chi\}}\big)\\
\leq&\mp(X_1^\eps\in\mathrm{supp}\chi)^\frac{1}{q}\|H_{(1,..,n)}(X_1^\eps,\chi(X_1^\eps))\|_p,
\end{align*}
where$\frac{1}{p}+\frac{1}{q}=1$. By Remark \ref{est H} we know that
$$\|H_{(1,...,n)}(X_1^\eps,\chi(X_1^\eps))\|_p\leq C_{p,q}\|\gamma_{X_1^\eps}^{-1}\|_\beta^m\|\bd X_1^\eps\|_{k,\gamma}^r\|\chi(X_1^\eps)\|_{k,q},$$
for some constants $\beta, \gamma>0$ and integers $k,m, r$. Thus, by Lemma \ref{th: Malliavin est subelliptic} we have
$$\lim_{\eps\downarrow0}\eps^2\log\|H_{(1,...,n)}(X_1^\eps,\chi(X_1^\eps))\|_p=0.$$

Finally by Theorem \ref{th: LDP}, a large deviation principle for $X_1^\eps$ ensures that for small $\eps$ we have
$$\mp(X_1^\eps\in\mathrm{supp}\chi)^\frac{1}{q}\leq e^{-\frac{1}{q\eps^2}(\inf_{y\in\mathrm{supp}\chi}d^2(y))}$$
which concludes the proof.  \hfill$\Box$

\noindent{\bf Sketch of the proof of (\ref{main claim 3})}.  Fix a point $y\in\mr^n$ and suppose that
$$\gamma:=\inf_{\Phi(h)=y, \det\gamma_{\Phi}(h)>0}\det\gamma_{\Phi}(h)>0.$$
Let  $\chi\in C_0^\infty(\mr^n), 0\leq\chi\leq1$ be a function such that $\chi$ is equal to one in a neighborhood of $y$,  and $g\in C^\infty(\mr), 0\leq g\leq1$, such that $g(u)=1$ if $|u|<\frac{1}{4}\gamma$, and $g(u)=0$ if $|u|>\frac{1}{2}\gamma$. Set $G_\eps=g(\det\gamma_{{X_1^\eps}})$. As before, we have
$$\me\chi(X_1^\eps)\delta_y(X_1^\eps)=\me G_\eps\chi(X_1^\eps)\delta_y(X_1^\eps)+\me(1-G_\eps)\chi(X_1^\eps)\delta_y(X_1^\eps).$$

First, it is easy to see that $\me G_\eps\chi(X_1^\eps)\delta_y(X_1^\eps)=0$ and  proceeding as in the proof of (\ref{main claim 2}) we obtain
\begin{align*}
\me(1-G_\eps)\chi(X_1^\eps)\delta_y(X_1^\eps)=&\me(\mathbf{1}_{\{X_1^\eps>y\}}H_{(1,...,n)}(X_1^\eps,(1-G_\eps)\chi(X_1^\eps)))\\
\leq&\me|H_{(1,...,n)}(X_1^\eps,(1-G_\eps)\chi(X_1^\eps))|\\
\leq&\me\big(|H_{(1,...,n)}(X_1^\eps\chi(X_1^\eps))|\mathbf{1}_{\{X_1^\eps\in\mathrm{supp}\chi, \det\gamma_{X_1^\eps}\geq\frac{1}{4}\gamma\}}\big)\\
\leq&\mp\left(X_1^\eps\in\mathrm{supp}\chi, \det\gamma_{X_1^\eps}\geq\frac{1}{4}\gamma\right)^{\frac{1}{q}}\|H_{(1,...,n)}(X_1^\eps,\chi(X_1^\eps))\|_p.
\end{align*}

Finally, by Lemma \ref{th: Malliavin est subelliptic}  and the large  deviation principle from Theorem \ref{th: LDP} for the couple $(X_1^\eps, \gamma_{X_1^\eps})$, we have for any $q>1$
\begin{align*}
\limsup_{\eps\downarrow0}\eps^2\log p_\eps(y)\leq&-\frac{1}{2q}\inf_{\Phi(h)\in\mathrm{supp}\chi, \det\gamma_{\Phi}(h)\geq\frac{1}{4}\gamma}\|h\|_\mathcal{H}^2\\
\leq&-\frac{1}{2q}\inf_{y\in\mathrm{supp}\chi}d^2_R(y).
\end{align*}
The proof is completed.\hfill$\Box$

\section{Small-time Kernel expansion}
\subsection{Laplace approximation}

Fix $H>\frac{1}{4}$ and consider equation (\ref{equ: SDE}). For the convenience of our discussion, in what follows, we write the above equation in the following form
\begin{align*}
X_t^\varepsilon=x+\varepsilon\int_0^t\sigma(X^\varepsilon_s)dB_s+\int_0^tb(\varepsilon,X_s^\varepsilon)ds,
\end{align*}
where $\sigma$ is a smooth $d\times d$ matrix and $b$ a smooth function from $\mr^+\times\mr^d$ to $\mr^d$. We also assume that $\sigma$ and $b$ have bounded derivatives to any order.

Fix $p>\frac{1}{H}$. Let $F$ and $f$ be two bounded infinitely Fr\'{e}chet differentiable  functionals on $\cC^{p-{\rm var}; [0,1]}([0,1], \R^d)$ with bounded derivatives (as linear operators) to any order. We are interested in studying the asymptotic behavior of 
$$J(\varepsilon)=\me\big[f(X^\varepsilon)\exp\{-F(X^\varepsilon)/\varepsilon^2\}\big],
\quad\quad \mathrm{as}\ \varepsilon \downarrow 0.$$ 

Recall for each $k\in \msh_H$,  $\Phi(k)$  is the deterministic It\^{o} map defined in (\ref{phi}). Set 
$$\Lambda(\phi)=\inf\{ \frac{1}{2}\|k\|_{\msh_H}, \phi=\Phi(k), k\in\msh_H\}.$$
 Throughout our discussion we make the following assumptions:
\begin{assumption}\label{assumption Laplace}
\ \\
\begin{itemize}
\item H 1: $F+\Lambda$ attains its minimum at finite number of paths $\phi_1, \phi_2, ..., \phi_n$ on $P(\mr^d)$.
\ \\
\item H 2: For each $ i\in\{1,2,...,n\}$, we have $\phi_i=\Phi(\gamma_i)$ and $\gamma_i$ is a non-degenerate minimum of the functional $F\circ\Phi+1/2\|\cdot\|^2_{\msh_H}$, i.e.:
$$ \forall k\in\msh_H\backslash\{0\},\quad d^2(F\circ\Phi+1/2\| \cdot\|^2_{\msh_H})(\gamma_i)k^2>0.$$
\end{itemize}
\end{assumption}

The following theorem is the main result of this section.

\begin{theorem}\label{Laplace}
Under the assumptions H 1 and H 2 above, we have
$$J(\varepsilon)=e^{-\frac{a}{\varepsilon^2}}e^{-\frac{c}{\varepsilon}}\bigg(\alpha_0+\alpha_1\varepsilon+...+\alpha_N\varepsilon^N+O(\varepsilon^{N+1})\bigg).$$
Here
$$a=\inf\{F+\Lambda(\phi), \phi\in P(\mr^d)\}=\inf\{F\circ\Phi(k)+1/2|k|^2_{\msh_H}, k\in\msh_H\}$$
and
$$c=\inf\big\{ dF(\phi_i)Y_i, i\in\{1,2,...,n\}\big\},$$
where $Y_i$ is the solution of
$$dY_i(s)=\partial_x\sigma(\phi_i(s))Y_i(s)d\gamma_i(s)+\partial_\varepsilon b(0,\phi_i(s))ds+\partial_x b(0,\phi_i(s))Y_i(s)ds$$
with $Y_i(0)=0$.
\end{theorem}

In what follows, we sketch the proof of the above Laplace approximation in the case $H>\frac{1}{2}$. Remarks on the rough case $\frac{1}{4}<H<\frac{1}{2}$ will be provided afterwards.

Without loss of generality, we may assume that $F+\Lambda$ attains its minimum at a unique path $\phi$.  There exists a $\gamma\in\msh_H$ such that $$\phi=\Phi(\gamma),\quad\quad \mathrm{and}\ \Lambda(\phi)=\frac{1}{2}\|\gamma\|^2_{\msh_H},$$ and
$$a\stackrel{\mathrm{def}}{=}\inf\{F+\Lambda(\phi), \phi\in P(\mr^d)\}=\inf\left\{F\circ\Phi(k)+\frac{1}{2}\|k\|^2_{\msh_H}, k\in\msh_H\right\}.$$
Moreover by assumption H 2,  for all non zero $k\in\msh_H$:
$$d^2(F\circ\Phi+\frac{1}{2}\|\ \|^2_{\msh_H})(\gamma)k^2>0.$$

Consider the following stochastic differential equation
$$Z^\varepsilon_t=x+\int_0^t\sigma(Z^\varepsilon_s)(\varepsilon dB_s+d\gamma_s)+\int_0^tb(\varepsilon,Z^\varepsilon_s)ds.$$
It is clear that $Z^0=\phi$. Denote $Z_t^{m,\varepsilon}=\partial_\varepsilon^mZ^\varepsilon_t$ and consider the Taylor expansion with respect to $\varepsilon$ near $\varepsilon=0$, we obtain
$$Z^\varepsilon=\phi+\sum_{j=0}^N\frac{g_j\varepsilon^j}{j!}+\varepsilon^{N+1}R^\varepsilon_{N+1},$$
where $g_j=Z^{j,0}.$ Explicitly, we have
\begin{align*}
dg_1(s)=\sigma(\phi_s)dB_s+\partial_x\sigma(\phi_s)g_1(s)d\gamma_s+\partial_xb(0,\phi_s)g_1(s)ds+\partial_\varepsilon b(0,\phi_s)ds.
\end{align*}

Now the proof is divided into the following steps.

\noindent{\underline {Step 1}}: By the large deviation principle, the sample paths that contribute to the asymptotics of $J(\varepsilon)$ lie in the neighborhoods of the minimizers  of $F+\Lambda$. More precisely, for $\rho>0$, denote by $B(\phi, \rho)$ the open ball (under $\lambda$-H\"{o}lder topology for a fixed $\lambda<H$) centered at $\phi$ with radius $\rho$. There exist $d>a$ and $\varepsilon_0>0$ such that for all $\varepsilon\leq \varepsilon_0$
$$\left|J(\varepsilon)-\me\left[f(X_T^\varepsilon)e^{-F(X_T^\varepsilon)/\varepsilon^2}, X^\varepsilon\in B(\phi,\rho)\right]\right|\leq e^{-d/\varepsilon^2}.$$
Hence, letting $$J_\rho(\varepsilon)=\me\left[f(X_T^\varepsilon)e^{-F(X_T^\varepsilon)/\varepsilon^2}, X^\varepsilon\in B(\phi,\rho)\right],$$
 to study the asymptotic behavior of $J(\varepsilon)$ as $\varepsilon\downarrow 0$, it suffices to study that of $J_\rho(\varepsilon)$.

\noindent{\underline {Step 2}}: Let $\theta(\varepsilon)=F(Z^\varepsilon)$ and write
$$\theta(\varepsilon)=\theta(0)+\varepsilon\theta'(0)+\frac{1}{2}\varepsilon^2\theta''(0)+\varepsilon^3R(\varepsilon).$$
By the Cameron-Martin theorem for fractional Brownian motions,  we have
\begin{align}\label{step 2}&J_\rho(\varepsilon)\\
=&\me\left\{f(Z^\varepsilon)\exp\left(-\frac{F(Z^\varepsilon)}{\varepsilon^2}\right)\exp\left(-\frac{1}{\varepsilon}\int_0^T\big((K^*_H)^{-1}(\dot{K_H^{-1}{\gamma}})\big)_sdB_s-\frac{\|\gamma\|^2_{\msh_H}}{2\varepsilon^2}\right); Z^\varepsilon\in B(\phi,\rho)\right\}\nonumber\\
=&\me\Bigg\{\exp\left[-\frac{1}{\varepsilon^2}\left(F(\phi)+\frac{1}{2}\|\gamma\|^2_{\msh_H}\right)\right]\exp\left[-\frac{\theta(0)'+\int_0^T\big((K^*_H)^{-1}(\dot{K_H^{-1}{\gamma}})\big)_sdB_s}{\varepsilon}\right]\nonumber\\
&\quad\quad\quad\exp\left[-\frac{1}{2}\theta''(0)\right]\cdot\left[ f(Z^\varepsilon)e^{-\varepsilon R(\varepsilon)}\right];Z^\varepsilon\in B(\phi,\rho)\Bigg\}.\nonumber
\end{align}

\noindent{\underline {Step 3}}:
It is clear that to prove Theorem \ref{Laplace}, it suffices to analyze the four terms in the expectation above. First of all, it is apparent that the first term ( of order -2) is
\begin{align}\label{e-2}\exp\left[-\frac{1}{\varepsilon^2}\left(F(\phi)+\frac{1}{2}\|\gamma\|^2_{\msh_H}\right)\right]=e^{{-\frac{a}{\varepsilon^2}}},\end{align}
which gives the leading term the Varadhan asymptotics.

The second term (of order -1) is deterministic. Indeed, since $\gamma$ is a critical point of $F\circ\Phi+1/2\|\cdot \|^2_{\msh_H}$ and note $\|k\|_{\msh_H}=\|K_H^{-1}k\|_\msh$, we have
\begin{align*}dF(\phi)(d\Phi(\gamma)k)=-\int_0^T\big((K^*_H)^{-1}\dot{(K_H^{-1}\gamma)}\big)_sdk_s.
\end{align*}
By the continuity of Young's integral with respect to the driving path, the above extends to
\begin{align*}dF(\phi)(d\Phi(\gamma)B)=-\int_0^T\big((K^*_H)^{-1}\dot{(K_H^{-1}\gamma)}\big)_sdB_s.
\end{align*}
On the other hand, note
$$\theta'(0)=dF(\phi)g_1,$$and
\begin{align*}g_1=d\Phi(\gamma)B+Y.\end{align*}
Here $Y$ is the solution of
$$dY_s=\partial_x\sigma(\phi_s)Y_sd\gamma_s+\partial_\varepsilon b(0,\phi_s)ds+\partial_x b(0,\phi_s)Y_sds,\quad\ Y(0)=0.$$
We obtain
\begin{align}\label{e-1}\exp\left[-\frac{\theta(0)'+\int_0^T\big((K^*_H)^{-1}(\dot{K_H^{-1}{\gamma}})\big)_sdB_s}{\varepsilon}\right]=\exp\left[-\frac{dF(\phi)Y}{\varepsilon}\right].\end{align}

For the third term (of order 0), one can show that there exists a $\beta>0$ such that
\begin{align}\label{e0}
\me \exp\left\{-(1+\beta)\left[\frac{1}{2}\theta''(0)\right]\right\}<\infty.
\end{align}
Let us emphasize that in order to show the above integrability of $\theta''(0)$, one needs to use assumption H2 and prove that $d^2F\circ\Phi(\gamma)(k^1,k^2)$ is Hilbert-Schmidt. For more details, we refer the reader to \cite{BO1}. Moreover, one can prove the following integrability of $R(\varepsilon)$.
\begin{lemma}\label{th: key lemma2}

There exist $\alpha>0$ and $\varepsilon_0>0$ such that
$$\sup_{0\leq \varepsilon\leq\varepsilon_0}\me\left(e^{(1+\alpha)|\varepsilon R(\varepsilon)|}; Z^\varepsilon\in B(\phi,\rho)\right)<\infty.$$\end{lemma}
Lemma \ref{th: key lemma2} and (\ref{e0}) allows us to analyze the third and forth terms and show 
\begin{align}\label{e higher}\me\big[f(Z^\varepsilon)e^{-\frac{1}{2}\theta''(0)-\varepsilon R(\varepsilon)}; Z^\varepsilon\in B(\phi,\rho)\big]=\sum_{m=0}^{N}\alpha_m\varepsilon^m+O(\varepsilon^{N+1}).\end{align} 
Finally, combining (\ref{step 2}), (\ref{e-2}), (\ref{e-1}), and (\ref{e higher}), the proof of Theorem \ref{Laplace} is complete. \hfill $\Box$

\begin{remark}\label{Rk Laplace frac} In application (see the next section), one may also be interested in an SDE which involves a fractional order term of $\varepsilon$,
\begin{align}\label{frac x}
X_t^\varepsilon=x+\varepsilon\int_0^t\sigma(X^\varepsilon_s)dB_s+\varepsilon^{\frac{1}{H}}\int_0^tb(\varepsilon,X_s^\varepsilon)ds. 
\end{align}
For this purpose, let us first introduce
\begin{align}\label{Inahama index}
\Lambda_1=\left\{n_1+\frac{n_2}{H}\big| n_1, n_2=0,1,2,...\right\},
\end{align}
the set of  fractional orders. Let $0=\kappa_0<\kappa_1<\kappa_2<\cdots$ be all elements of $\Lambda_1$ in increasing order. When $H>\frac{1}{2}$, we have
\begin{align}\label{kappa}(\kappa_0, \kappa_1,\kappa_2,\kappa_3,\kappa_4,...)=(0, 1, \frac{1}{H}, 2,1+\frac{1}{H},...).\end{align}
Set 
$$\Lambda_2=\{\kappa-2|\kappa\in\Lambda_1\backslash\{0\}\},$$
and define
$$\Lambda_3=\{a_1+a_2+\dots+a_m|m\in\mn_+\ \mathrm{and}\ a_1,...,a_m\in\Lambda_1\}$$
and
$$\Lambda_3'=\{a_1+a_2+\dots+a_m|m\in\mn_+\ \mathrm{and}\ a_1,...,a_m\in\Lambda_2\}.$$
Finally let
$$\Lambda_4=\{a+b| a\in\Lambda_3, b\in\Lambda_3'\}$$
and denote by $\{0=\lambda_0<\lambda_1<\lambda_2<\dots\}$ all the elements of $\Lambda_4$ in increasing order. Let us note that  the set $\Lambda_3$ characterizes the powers of $\varepsilon$ coming from the term $f(Z^\varepsilon)$ in (\ref{step 2}) and $\Lambda_3'$ characterizes that of $e^{-\varepsilon R(\varepsilon)}$.

Similar as before, we consider
\begin{align}\label{frac z}Z^\varepsilon_t=x+\int_0^t\sigma(Z^\varepsilon_s)(\varepsilon dB_s+d\gamma_s)+\varepsilon^{\frac{1}{H}}\int_0^tb(\varepsilon,Z^\varepsilon_s)ds.\end{align}
It can be proved that $Z^\epsilon$ has the following expansion in $\varepsilon$,
$$Z^\varepsilon=\phi+\sum_{j=0}^Ng_{\kappa_j}\varepsilon^{\kappa_j}+\varepsilon^{\kappa_N+1}R^\varepsilon_{\kappa_N+1}.$$
Note that in (\ref{kappa}), indices up to degree two are $(0, 1, 1/H, 2)$. There is an extra term $1/H$ compared to the case without fractional order. Hence when plugging  (\ref{frac z}) into Step 2 of the proof of Theorem \ref{Laplace}, there is an extra (but deterministic) term
$$\exp\left\{-\frac{dF(\Phi)g_{\kappa2}}{\varepsilon^{2-\frac{1}{H}}}\right\},$$
where $g_{\kappa_2}$ satisfies
$$d{g_{\kappa_2}}(s)=\partial_x\sigma(\phi_s){g_{\kappa_2}}(s)d\gamma_s+b(0,\phi_s)ds,\quad g_{\kappa_2}(0)=0.$$
It is not hard to see that the other terms up to degree two remain the same, and that although higher order terms are different they could be handled similarly as before. Hence we obtain
\begin{theorem}\label{Laplace 2} Let $X^\varepsilon$ satisfy (\ref{frac x}). we have
\begin{align*}\me\big[f(X^\varepsilon)e^{-F(X^\varepsilon)/\varepsilon^2}\big]
=e^{-\frac{a}{\varepsilon^2}}e^{-\frac{c}{\varepsilon}}\exp\left\{-\frac{d}{\varepsilon^{2-\frac{1}{H}}}\right\}\bigg(\alpha_{\lambda_0}+\alpha_{\lambda_1}\varepsilon^{\lambda_1}+...+\alpha_{\lambda_N}\varepsilon^{\lambda_N}+O(\varepsilon^{\lambda_{N+1}})\bigg).\end{align*}
Here
\begin{align*} &a=\inf\{F\circ\Phi(k)+1/2|k|^2_{\msh_H}, k\in\msh_H\},\\ &c=dF(\phi)Y,\quad \mathrm{and}\quad d=dF(\phi)g_{\kappa_2},\end{align*}
where $Y$ and $g_{\kappa_2}$ satisfiy
$$dY(s)=\partial_x\sigma(\phi_i(s))Y(s)d\gamma(s)+\partial_\varepsilon b(0,\phi(s))ds+\partial_x b(0,\phi(s))Y(s)ds,\quad Y(0)=0,$$
and
$$d{g_{\kappa_2}}(s)=\partial_x\sigma(\phi_s){g_{\kappa_2}}(s)d\gamma_s+b(0,\phi_s)ds,\quad g_{\kappa_2}(0)=0.$$
\end{theorem}
\end{remark}

\begin{remark}
Theorem \ref{Laplace} for the rough case $\frac{1}{4}<H<\frac{1}{2}$ was proved by Inahama \cite{Inahama1}. In this case, equation is understood in the rough path sense. Thanks to Proposition \ref{prop:imbed-bar-H}, equations for $g_i$ and $R_i$ are understood as Young's paring.

In \cite{Inahama1} the author also discussed RDEs with fractional orders of $\varepsilon$, in which the index set $\Lambda_1$ was introduced. The main idea of the proof for the rough case is the same as that outlined above. But the major difficulty is to show that $d^2F\circ\Phi(\gamma)(k^1,k^2)$ is Hilbert-Schmidt. This is easier when $H>\frac{1}{2}$, since in this case $\partial_tK(t,s)$ is integrable,  and one can easily obtain a nice representation for $d^2F\circ\Phi(\gamma)(k^1,k^2)$.
\end{remark}

\subsection{Expansion of the density function} Consider
\begin{align}\label{SDE with drift}
X_t=x+\sum_{i=1}^d\int_0^tV_i(X_s)dB^i_s+\int_0^tV_0(X_s)ds.
\end{align}
We are interested in studying the small-time asymptotic behavior of $X_t$.  It is clear that by the self-similarity of $B$,  this is equivalent to studying the asymptotic behavior of $X_1^\varepsilon$ (for small $\varepsilon$) which satisfies
$$X_t^\varepsilon=x+\sum_{i=1}^d\varepsilon\int_0^tV_i(X_s)dB^i_s+\varepsilon^{\frac{1}{H}}\int_0^tV_0(X_s)ds.$$

In what follows, we use the Laplace approximation to obtain a short time asymptotic expansion for the density of $X_1^\varepsilon$ in the case when $H>\frac{1}{2}$. For this purpose, we need the following assumption.
\begin{assumption}\label{H}

\
\begin{itemize}
\item A 1: For every $x \in \mathbb{R}^d$, the vectors $V_1(x),\cdots,V_d(x)$ form a basis of $\mathbb{R}^d$.
\item A 2: There exist smooth and bounded functions $\omega_{ij}^l$ such that:
\[
[V_i,V_j]=\sum_{l=1}^d \omega_{ij}^l V_l,
\]
and
\[
\omega_{ij}^l =-\omega_{il}^j.
\]
\end{itemize}
\end{assumption}

Assumption A1 is the standard ellipticity condition. Due to the second assumption A2, the geodesics are easily described. If $k: \mathbb{R}_{\ge 0} \rightarrow \mathbb{R}$ is a $\alpha$-H\"older  path with $\alpha >1/2$ such that $k(0)=0$, we denote by $\Phi(x,k)$ the solution of the ordinary differential equation:
\[
x_t=x+\sum_{i=1}^d \int_0^t V_i(x_s) dk^i_s.
\]
Whenever there is no confusion, we always suppress the starting point $x$ and denote it simply by $\Phi(k)$ as before. Then we have (see Lemma 4.2 in \cite{BO1})
\begin{lemma}
$\Phi(x,k)$ is a geodesic if and only if $k(t)=t u$ for some $u \in \mathbb{R}^d$.
\end{lemma}
As a consequence of the previous  lemma, we then have the following key result (Proposition 4.3 in \cite{BO1}):

\begin{proposition}\label{th: ministance}
Let $T>0$. For $x,y \in \mathbb{R}^d$,
\[
\inf_{k \in \msh_H, \Phi_T(x,k)=y } \| k \|^2_{\msh_H}=\frac{d^2(x,y)}{T^{2H}} . 
\]
\end{proposition}

\begin{lemma}\label{Func F}
For any $x\in\mr^d$,  there exists a neighborhood $V$ of $x$ and a bounded smooth function $F(x,y,z)$ on $V\times V\times \mr^d$ such that:

(1) For any $(x,y)\in V\times V$ the infimum 
$$\inf\left\{F(x,y,z)+\frac{d(x,z)^2}{2}, z\in M\right\}=0$$
is attained at the unique point $y$. Moreover, it is a non-degenerate minimum. Hence there exists a unique $k^0\in\msh_H$ such that
(a): $\Phi_1(x_0,k^0)=y_0$;
(b): $d(x_0,y_0)=\|k^0\|_{\msh_H}$; and (c): $k^0$ is a non-degenerate minimum of the functional: $k\rightarrow F(\Phi_1(x_0,k))+1/2\|k\|^2_{\msh_H}$ on $\msh_H$.

(2) For each $(x,y)\in V\times V$, there exists a ball centered at $y$ with radius $r$ independent of $x,y$ such that $F(x,y,\cdot)$ is a constant outside of the ball.
\end{lemma}

Let $F$ be in the above lemma and $p_\varepsilon(x,y)$ the density function of $X^\varepsilon_1$.  By the inversion of Fourier transformation we have
\begin{align}\label{Fourier}
p_\varepsilon(x,y)e^{-\frac{F(x,y,y)}{\varepsilon^2}}&=\frac{1}{(2\pi)^d}\int e^{-i\zeta\cdot y}d\zeta\int e^{i\zeta\cdot z}e^{-\frac{F(x,y,z)}{\varepsilon^2}}p_\varepsilon (x,z)dz\\
&=\frac{1}{(2\pi\varepsilon)^d}\int e^{-i\frac{\zeta\cdot y}{\varepsilon}}d\zeta\int e^{i\frac{\zeta\cdot z}{\varepsilon}}e^{-\frac{F(x,y,z)}{\varepsilon^2}}p_\varepsilon(x,z)dz\nonumber\\
&=\frac{1}{(2\pi\varepsilon)^d}\int d\zeta \me_x\left(e^{\frac{i\zeta\cdot(X_1^\varepsilon-y)}{\varepsilon}}e^{-\frac{F(x,y,X^\varepsilon_1)}{\varepsilon^2}}\right).\nonumber
\end{align}
It is clear that by  applying Laplace approximation to the expectation in the last equation above and  
switching the order of integration (with respect to $\zeta$) and summation, we obtain an asymptotic expansion for the the density function $p_\varepsilon(x,y)$.  

\begin{remark} One might wonder why not constructing, for each fixed $x,y$, a function $F$ which minimizes (at $z=y$)
$$F(x,y,z)+\frac{D(x,z)^2}{2}$$
in Lemma \ref{Func F}, where 
$$D^2(x,y)=\inf_{k \in \msh_H, \Phi_1(x,k)=y } \| k \|^2_{\msh_H}. $$
After all $D(x,y)$ seems the natural ``distance'' for the system (\ref{SDE with drift}), instead of the Riemannian distance $d(x,y)$. The problem with $D(x,y)$ is that it is not clear weather it is differentiable, while the construction of $F$ in Lemma \ref{Func F} needs some differentiability of $D(x,y)$. This is indeed one of the reasons why we impose the structure assumption A2 so that $D(x,y)=d(x,y)$ (content of Proposition \ref{th: ministance}). With this identification, we know $D(x,y)$ is smooth for all $x\not=y$. 
\end{remark}

\begin{remark}
In order to show Proposition \ref{th: ministance}, we used the fact that $\partial K(t,s)/\partial t$ is integrable, which is only true for the smooth case $H>\frac{1}{2}$. Hence although Inahama proved the Laplace approximation for $\frac{1}{4}<H<\frac{1}{2}$ in \cite{Inahama1}, we can not repeat the proof in this section to produce an expansion of the density function for the rough case.
\end{remark}

Recall the definition of $\Lambda_1$ in Remark \ref{Rk Laplace frac} and similarly set 
$$\Lambda_2=\{\kappa-1|\kappa\in\Lambda_1\backslash\{0\}\}$$
and 
$$\Lambda_2'=\{\kappa-2|\kappa\in\Lambda_1\backslash\{0\}\}.$$
Next define
$$\Lambda_3=\{a_1+a_2+\dots+a_m|m\in\mn_+\ \mathrm{and}\ a_1,...,a_m\in\Lambda_2\}.$$
and
$$\Lambda_3'=\{a_1+a_2+\dots+a_m|m\in\mn_+\ \mathrm{and}\ a_1,...,a_m\in\Lambda_2'\}.$$
Finally, set
$$\Lambda_4=\{a+b| a\in\Lambda_3, b\in\Lambda_3'\}$$
and denote by $\{0=\lambda_0<\lambda_1<\lambda_2<\dots\}$ all the elements of $\Lambda_4$ in increasing order. Similar as before, powers of $\epsilon$ in the index set $\Lambda_3$ comes from the term $\exp\left\{{i\zeta\cdot(X_1^\varepsilon-y)}/{\varepsilon}\right\}$ in (\ref{Fourier}) and powers in $\Lambda_3'$ comes from $\exp\{-{F(x,y,X^\varepsilon_1)}/{\varepsilon^2}\}$.

Our main result of this section is the following (by letting $\varepsilon=t^H$).
\begin{theorem}
Fix $x\in \mr^d$. Suppose the Assumption \ref{H} is satisfied, then in a neighborhood $V$ of $x$, the density function $p(t;x,y)$ of $X_t$ in (\ref{SDE with drift}) has the following asymptotic expansion near $t=0$
\begin{align*}
p(t;x,y)=\frac{1}{(t^H)^{d}}e^{-\frac{d^2(x,y)}{2t^{2H}}+\frac{\beta}{t^{2H-1}}}\bigg(\sum_{i=0}^N c_i(x,y)t^{\lambda_iH}+r_{N+1}(t,x,y)t^{\lambda_{N+1}H}\bigg),\quad\quad y\in V.
\end{align*}
Here $\beta$ is some constant, $d(x,y)$ is the Riemannian distance between $x$ and $y$ determined by $V_1,...,V_d$. Moreover, we can chose $V$ such that $c_i(x,y)$ are $C^\infty$ in  $ V\times V\subset \mr^d\times \mr^d$,   and for all multi-indices $\alpha$ and $\beta$ 
$$\sup_{t\leq t_0}\sup_{(x,y)\in V\times V}|\partial^\alpha_x\partial^\beta_y r_{N+1}(t,x,y)|<\infty$$
for some $t_0>0$.
\end{theorem}

\begin{remark}
Differentiability of $c_i(x,y), r_{N+1}$ in the above theorem and legitimacy of Fourier inversion in (\ref{Fourier}) is obtained by Malliavin calculus and some uniform estimates of the coefficients in the Laplace approximation. We refer the reader to \cite{BO1} for details.
\end{remark}

\begin{remark}
Our result assumes the  ellipticity condition and a strong structure condition (Assumption \ref{H}). Later Inahama \cite{Inahama2} proved the kernel expansion under some mild conditions on the vector fields (also in the smooth case $H>\frac{1}{2}$). He takes a different approach and uses Watanabe distribution theory. On the other hand, the smoothness of coefficient and the uniform estimate for the remainder terms in the expansion are not provided in \cite{Inahama2}.
\end{remark}

\section{Application to mathematical finance}
Fractional Brownian motions has been used in financial models to introduce memory. In this section, we give two examples of such models and  remark on how the methods and results in the previous sections could be applied to the study of such models.

\subsection{One dimensional models} Memories can be introduced to stock price process directly. In particular, the so-called fractional Black and Scholes model is given by 
\begin{align}\label{fBS}
S_t=S_0\exp\left(\mu t+\sigma B_t^H-\frac{\sigma^2}{2}t^{2H}\right),
\end{align}
where $B^H$ is a fractional Brownian motion with Hurst parameter $H$, $\mu$ the mean rate of return and $\sigma>0$ the volatility.  Let $r$ be the interest rate. The price for the risk-free  bond is given by $e^{rt}$.

More generally, one can also consider a fractional local volatility model 
\begin{align*}
dS_t=S_t(\mu dt+\sigma(S_t) dB_t^H).
\end{align*}
Here the stochastic integration with respect to $B^H$ could be understood in the sense of rough path theory. After a simple change of variable $X_t=\log S_t$, one obtains
$$dX_t=\mu dt+\sigma(e^{X_t})dB_t^H.$$
There has been an intensive study recently of option prices and implied volatilities for options with short maturity (e.g. \cite{Busca2:2004}, \cite{Gatheral:2009}, \cite{Forde09}, \cite{FengFordeFouque}). Since the above equation is a special case of (\ref{SDE with drift}), we can use the results obtained in the previous sections to obtain short-time asymptotic behavior of such models.

A drawback of the finance models discussed above is that they lead to the existence of arbitrage opportunities. For example, let the couple $(\alpha_t, \beta_t), t\in[0,T]$ be a  portfolio with $\alpha_t$ the amount of bonds and $\beta_t$ the amount of stocks at time $t$. When $H>\frac{1}{2}$, one can construct an arbitrage in the fractional Black and Scholes model by (for simplicity, we assume $\mu=r=0$)
\begin{align*}
\beta_t=S_t-S_0,\quad \mathrm{and}\ \ \alpha_t=\int_0^t\beta_tdS_t-\beta_tS_t.
\end{align*}

\subsection{Stochastic volatility models} Stochastic volatility models was introduced to capture both the volatility smile and the correct dynamics of the volatility smile (see \cite{Hagan} for instance). For these models, molding the volatility process is one of the key factors. In \cite{Comte}, the authors proposed a long memory specification of the volatility process in order to capture the steepness of long term volatility smiles without over increasing the short run persistence. 

The following stochastic volatility model based on the fractional Ornstein-Uhlenbeck process provides another way introducing long memory to the volatility process:
$$dS_t=\mu S_tdt+\sigma_tS_tdW_t,$$
where $\sigma_t=f(Y_t)$ and $Y_t$ is a fractional Ornstein-Uhlenbeck process:
$$dY_t=\alpha(m-Y_t)dt+\beta_tdB_t^H.$$
In the above $W_t$ is a standard Brownian motion and $B^H_t$ an independent (of $W_t$)  fractional Brownian motion with Hurst parameter $H>\frac{1}{2}$. Examples of functions $f$ are $f(x)=e^x$ and $f(x)=|x|$.

Comte and Renault \cite{CR} studied this type of stochastic volatility models which introduces long memory and mean reverting in the Hall and White setting \cite{HW}. The long memory property allows this model to capture the well-documented evidence of persistence of the stochastic feature of Black and Scholes implied volatilities when time to maturity increases.

Unlike one dimensional models mentioned above, the fractional Ornstein-Uhlenbeck model is arbitrage free since the stock price process is driven by a standard Brownian motion. In \cite{Hu}, Hu has proved that for this model, market is incomplete and the martingale measures are not unique. If we set $\gamma_t=(r-\mu)/\sigma_t$ and
$$\frac{d\mq}{d\mp}=\exp\left(\int_0^T\gamma_tdW_t-\frac{1}{2}\int_0^T|\gamma_t|^2dt\right).$$
Then $\mq$ is the minimal martingale measure associated with $\mp$. Moreover, the risk minimizing-hedging price at $t=0$ of an European call option with payoff $(S_T-K)^+$ is given by 
$$C_0=e^{-rT}\me_{\mq}(S_T-K)^+.$$

The fractional Ornstein-Uhlenbeck model takes a generalized form of equation (\ref{SDE with drift}) that is studied in the previous sections. It is a system of SDEs driven by fractional Brownian motions, but with varying Hurst parameter $H$. We believe that the methods discussed above can be extended to study small-time asymptotics of these models.

\end{document}